\documentclass[a4paper,12pt]{article}     
\usepackage{amsmath,amscd,amssymb}        
\usepackage{latexsym}                     
\usepackage[english, german]{babel}                

\usepackage{theorem}                 

\newtheorem{theorem}{Theorem}
\newtheorem{corollary}{Corollary}
\newtheorem{proposition}{Proposition}

\newenvironment{example}
{\smallskip\noindent{\bf Example\/}.}{\smallskip\par}

\newenvironment{remark}
{\smallskip\noindent{\bf Remark\/}.}{\smallskip\par}

\newenvironment{proof}{\begin{ProofwCaption}{Proof}}{\end{ProofwCaption}}
\newenvironment{proof*}[1]{\begin{ProofwCaption}{{#1}}}{\end{ProofwCaption}}
\newenvironment{ProofwCaption}[1]%
  {\addvspace\theorempreskipamount \noindent{\it #1.}\rm}%
  {\qed \par \addvspace\theorempostskipamount}
\newcommand{\qedsymbol}{\mbox{$\Box$}}
\newcommand{\qed}{\hfill\qedsymbol}

\newcommand{\CC}{{\mathbb C}}

\newcommand{\PP}{{\mathbb P}}

\newcommand{\ZZ}{{\mathbb Z}}

\newcommand{\calO}{{\cal O}}

\newcommand{\calJ}{{\cal J}}

\newcommand{\calE}{{\cal E}}

\newcommand{\calX}{{\cal X}}

\newcommand{\calD}{{\cal D}}

\newcommand{\eps}{\varepsilon}

\newcommand{\vv}{\underline{v}}
\newcommand{\ww}{\underline{w}}

\newcommand{\ttt}{\underline{t}}
\newcommand{\mm}{\underline{m}}

\title{On divisorial filtrations associated with Newton diagrams}
\author{W.~Ebeling and S.~M.~Gusein-Zade
\thanks{Partially supported by the DFG (Eb 102/7--1), RFBR--10-01-00678,
NSh--8462.2010.1.
AMS 2010 Math. Subject Classification: 32S05, 14M25, 16W70.
}
}
\date{}

\begin{document}
\selectlanguage{english}

\maketitle

\begin{abstract}
We consider divisorial filtration on the rings of functions on hypersurface singularities associated with Newton diagrams and their analogues for plane curve singularities. We compute the multi-variable Poincar\'e series for the latter ones.
\end{abstract}

\section*{Introduction} 
A multi-index filtration on the ring $\calO_{V,0}=\calO_{\CC^n,0}/(f)$ of functions on a hypersurface singularity $(V,0)=\{f=0\}$ defined by the Newton diagram $\Gamma=\Gamma_f$ of the germ $f$ was considered in \cite{JS}. The initial idea was to look for a filtration corresponding to a Newton diagram for which the Poincar\'e series could be computed and compared with the corresponding monodromy zeta function. This was inspired by the coincidence of Poincar\'e series and monodromy zeta functions in some cases (e.g.\ in \cite{IJM}) and relations between them in some other cases (e.g.\ in \cite{MRL}).
A somewhat natural filtration on the ring $\calO_{V,0}$ corresponding to the Newton diagram $\Gamma=\Gamma_f$ is the divisorial filtration defined by the divisors in a toric resolution of $f$ corresponding to the facets of the diagram. However, at that moment the divisorial valuation was regarded as being complicated to treat. The filtration defined in \cite{JS} was regarded as a certain ``simplification'' of the divisorial one. In fact this seems to be not the case. It is rather complicated to compute the Poincar\'e series of that filtration and moreover the assertion of Theorem~1 of \cite{JS} for $s>2$ appeared to be wrong. Another filtration corresponding to a Newton diagram was considered in \cite{L}.

Here we discuss an analogue of the divisorial valuation corresponding to a Newton diagram, describe its generalization for plane curve singularities, and compute the Poincar\'e series of the latter one.

For a germ $(V,0)$ of a complex analytic variety, let $\pi:(\calX, \calD)\to(V,0)$ be a resolution of $(V,0)$ with the exceptional divisor $\calD=\pi^{-1}(0)$ being a normal crossing divisor on $\calX$.
For an irreducible component $\calE$ of $\calD$ and for $g\in{\calO}_{V,0}$, let $v_{\calE}(g)$ be the order of the zero of the lifting $\widetilde{g}=g\circ\pi$ of the germ $g$ to the space $\calX$ of the resolution along $\calE$. 
The function $v_{\calE}:\calO_{V,0}\to\ZZ_{\ge 0}\cup\{+\infty\}$ is called a {\em divisorial valuation} on $\calO_{V,0}$. 
One can consider the multi-index filtration defined by a collection $\calE_1$, \dots, $\calE_r$ of components of the exceptional divisor:
\begin{equation}\label{filt}
J(\vv)=\{g\in\calO_{V,0}: \vv(g)\ge\vv\}\,,
\end{equation}
where $\vv=(v_1, \ldots, v_r)\in \ZZ_{\ge 0}^r$,
$\vv(g)=(v_1(g), \ldots, v_r(g))$,
$v_i(g)=v_{\calE_i}(g)$,
$\vv'=(v_1', \ldots, v_r')\ge\vv$ if and only if
$v_i'\ge v_i$ for $i=1, \ldots, r$. 
This filtration is called a {\em divisorial} one. 
The notion of the Poincar\'e series of a multi-index filtration was introduced in \cite{CDK} (see also \cite{IJM}).
In \cite{IJM} it was explained that the Poincar\'e series of a filtration defined by a formula like (\ref{filt}) is equal to the integral with respect to the Euler characteristic
\begin{equation}\label{Poincare}
P_{\{v_i\}}(\ttt)=\int_{\PP\calO_{V,0}}\ttt^{\vv(g)}d\chi
\end{equation}
over the projectivization $\PP\calO_{V,0}$ of $\calO_{V,0}$ ($\ttt=(t_1, \ldots , t_r)$, $\ttt^{\vv}=t_1^{v_1} \cdots t_r^{v_r}$). In this integral, $t_i^\infty$ has to be assumed to be equal to zero.
Also in \cite{IJM} it was shown that the Poincar\'e series of the divisorial filtration corresponding to all the components of the exceptional divisor of a resolution (uniformization) of a plane curve singularity $(C,0)=\{f=0\}\subset(\CC^2,0)$ (that is to all the components of the curve $(C,0)$) coincides with the Alexander polynomial (in several variables) of the corresponding link $C\cap S^3_{\eps}\subset S^3_{\eps}$, where $S^3_{\eps}$ is the sphere of small radius $\eps$ centred at the origin in $\CC^2$. 
(The Alexander polynomial becomes the monodromy zeta function of the left hand side $f$ of the equation of the curve $(C,0)$ after identification of all the variables.)

For the definition of a multi-index filtration by the formula (\ref{filt}), it is not necessary to assume that all the $v_i:\calO_{V,0}\to\ZZ_{\ge 0}\cup\{+\infty\}$ are valuations (i.e. that they satisfy the condition $v_i(g_1g_2)=v_i(g_1)+v_i(g_2)$). It is sufficient to require that all of them are so called {\em order functions}. This means that they satisfy the condition
$v_i(g_1+g_2)\ge \min{\{v_i(g_1),v_i(g_2)\}}$,
but, in general not the condition $v_i(g_1g_2)=v_i(g_1)+v_i(g_2)$. We shall use order functions to define the filtrations below.

\section{Divisorial filtration corresponding to a Newton diagram}
Let $f\in \calO_{\CC^n,0}$ be a function germ non-degenerate with respect to its Newton diagram $\Gamma=\Gamma_f$.
Let $p:(X,D)\to (\CC^n,0)$ be a toric resolution of $f$ corresponding to the Newton diagram $\Gamma$.
The facets of $\Gamma$ correspond to some components (say, $E_1$, \dots, $E_r$) of the exceptional divisor $D$. Let $(V,0)=\{f=0\}$, let $\widetilde{V}$ be the strict transform of the hypersurface singularity $V$, and let $\calE_i:=\widetilde{V}\cap E_i$.

For $n\ge 3$ the $\calE_i$ are the irreducible components of the exceptional divisor $\calD=D\cap\widetilde{V}$ of the resolution $p_{\vert\widetilde{V}}:(\widetilde{V}, \calD)\to (V,0)$. Thus one can consider the divisorial valuations $v_i$ defined by these components and the corresponding (multi-index) filtration on $\calO_{V,0}$. For $n=2$ the intersections $\calE_i$ are not, in general, irreducible (if the corresponding facets of $\Gamma$ have integer points in their interiors). Therefore for $n=2$ the corresponding definition has to be modified.

Let us first reformulate the definition of the divisorial valuations (for $n\ge 3$) in terms of the Newton diagram $\Gamma$. Let $\gamma_1$, \dots, $\gamma_r$ be the facets of the diagram $\Gamma$ and let $\ell_i(\bar{k})=c_i$ be the reduced equation of the facet $\gamma_i$, $i=1, \ldots, r$. This means that $\ell_i(\bar{k})=a_{i1}k_1+\ldots+a_{in}k_n$ ($\bar{k}=(k_1, \ldots, k_n)$), where $a_{i1}$, \dots, $a_{in}$ are positive integers with greatest common divisor equal to $1$.

For $g\in\calO_{\CC^n,0}[x_1^{-1}, \ldots, x_n^{-1}]$, $g(\bar{x})=\sum\limits_{\bar{k}}c_{\bar{k}}\bar{x}^{\bar{k}}$ ($\bar{x}=(x_1, \ldots , x_n)$), let $u_i(g):=\min\limits_{\bar{k}: c_{\bar{k}}\ne 0}\ell_i(\bar{k})$, and let $g_{\gamma_i}(\bar{x})=\sum\limits_{\bar{k}: \ell_i(\bar{k})=u_i(g)}c_{\bar{k}}\bar{x}^{\bar{k}}$. For $g\in \calO_{\CC^n,0}/(f)$ (or rather for $g\in \calO_{\CC^n,0}$) let us define $\widehat{v}_i(g)$ by
\begin{equation}\label{v_hat}
\widehat{v}_i(g)=\sup_{h\in \calO_{\CC^n,0}[x_1^{-1}, \ldots, x_n^{-1}]} u_i(g+hf)\,.
\end{equation}

One can see that, for $n=2$, $\widehat{v}_i:\calO_{\CC^n,0}/(f)\to \ZZ_{\ge 0}\cup\{+\infty\}$ is not, in general, a valuation, but only an order function.

\begin{example}
Let $f(x,y)=y^3+y^2x-x^5$ and let $\gamma_1$ be the facet of $\Gamma_f$ defined by the equation $2k_y+k_x=5$. Let $g_1(x,y)=y+x^2$, $g_2(x,y)=y-x^2$. One has $\widehat{v}_1(g_i)=u_1(g_i)=2$ for $i=1,2$, but $\widehat{v}_1(g_1g_2) = u_1(g_1g_2-x^{-1}f) = u_1(-y^3x^{-1})=5$.
\end{example}

\begin{remark}
One can see that this definition resembles the definition used in \cite{JS} where similar order functions were defined by equation (\ref{v_hat}) with $\calO_{\CC^n,0}[x_1^{-1}, \ldots, x_n^{-1}]$ substituted by $\calO_{\CC^n,0}$.
\end{remark}

\begin{proposition}\label{prop1}
For $n\ge 3$, $i=1, \ldots , r$, and $g\in \calO_{\CC^n,0}$ one has
$$
\widehat{v}_i(g)={v}_i(g)\,.
$$
\end{proposition}

\begin{proof}
The claim follows from the following statements:
\begin{enumerate}
\item[1)] $v_i(g)\ge \widehat{v}_i(g)$;
\item[2)] if $f_{\gamma_i} \!\! \not\vert \, g_{\gamma_i}$, then $v_i(g)=u_i(g)$;
\item[3)] if $f_{\gamma_i}\vert g_{\gamma_i}$, then there exists $h\in \calO_{\CC^n,0}[x_1^{-1}, \ldots, x_n^{-1}]$ such that $u_i(g+hf)>u_i(g)$.
\end{enumerate}

Indeed, by iterated applications of 2) and 3) one obtains that either $\widehat{v}_i(g)=\infty$ or there exists $h\in \calO_{\CC^n,0}[x_1^{-1}, \ldots, x_n^{-1}]$ such that $v_i(g)=u_i(g+hf)$. Therefore 
$\widehat{v}_i(g)=\sup\limits_{h\in \calO_{\CC^n,0}[x_1^{-1}, \ldots, x_n^{-1}]}u_i(g+hf)\ge v_i(g)$ and 1) implies the assertion.

Statement 1) follows from the facts that: $u_i(g)$ is the order of vanishing of the lifting $g\circ\pi$ of $g$ along $E_i$; $v_i(g)$ is the order of vanishing of $g\circ\pi_{\vert\widetilde{V}}$ along $\calE_i\subset E_i$ and therefore $v_i(g)\ge u_i(g)$; $v_i(g)=v_i(g+hf)$ for any $h\in \calO_{\CC^n,0}[x_1^{-1}, \ldots, x_n^{-1}]$.

If $f_{\gamma_i} \!\! \not\vert \, g_{\gamma_i}$, then the intersection $\widetilde{\{g=0\}}\cap E_i$ of the strict transform $\widetilde{\{g=0\}}$ with the component $E_i$ does not contain $\calE_i$. Therefore the order of vanishing of $g\circ\pi_{\vert\widetilde{V}}$ along $\calE_i$ coincides with the order of vanishing of $g\circ\pi$ along $E_i$, equal to $u_i(g)$. This gives 2).

If $g_{\gamma_i}=hf_{\gamma_i}$ ($h\in \calO_{\CC^n,0}[x_1^{-1}, \ldots, x_n^{-1}]$), then $(g-hf)_{\gamma_i}$ contains with non-zero coefficients only monomials ${\bar{x}}^{\bar{k}}$ with $\ell_i(\bar{k})> u_i(\bar{k})$. This gives 3).
\end{proof}

As it was mentioned above, for $n=2$ the intersections $\calE_i=E_i\cap \widetilde{V}$ may be reducible: i.e. consist of several points. In this case there is no divisorial valuation associated to $\calE_i$. Let us modify (generalize) the definition of a divisorial valuation in the following way. Let $\calE = \bigcup\limits_{j=1}^s \calE^{(j)}$ be the union of irreducible components of the exceptional divisor $\calD$ of the resolution $\pi:(\widetilde{V}, \calD)\to(V,0)$ and for 
$g\in \calO_{\CC^n,0}$ let 
$$
v_{\calE}(g):=\min\limits_{j=1,\ldots, s}v_{\calE^{(j)}}(g)\,.
$$
The function $v_{\calE}:\calO_{V,0}\to\ZZ_{\ge 0}\cup\{+\infty\}$ is not, in general, a valuation (for $s>1$), but an order function. The number $v_{\calE}(g)$ can also be defined as the minimum over all arcs $\gamma$ on $\widetilde{V}$ at points of $\calE$ of the order of $g$ along $\gamma$.

One can easily see that this definition gives order functions $v_i$ on $\calO_{V,0}$ corresponding to the facets of the Newton diagram $\Gamma=\Gamma_f$ for $n=2$ as well so that Proposition~\ref{prop1} also holds in this case.

\section{Plane curve singularities}
Here we consider analogues 
of the order functions $v_i$ corresponding to the facets of the Newton diagram $\Gamma$ (for $n=2$) for plane curve singularities not associated with Newton diagrams (say, for those whose components may have more than one Puiseux pair). We compute the Poincar\'e series of the corresponding filtration and give its specialization for the filtration defined by a Newton diagram. It seem to be less involved to carry out computations in this way than to produce them directly by considering Newton diagrams.

Let $(C,0)\subset(\CC^2,0)$ be a plane curve singularity with an embedded resolution $\pi:(X,D)\to(\CC^2,0)$ such that
\begin{enumerate}
 \item[1)] $C$ is the union of irreducible components $C=\bigcup\limits_{i,j}C_{ij}$, where $i=1, \ldots, r$, $j=1, \ldots, s_i$ ($s_i>0$);
 \item[2)] for each $i$ the strict transforms $\widetilde{C}_{i1}$, \dots $\widetilde{C}_{is_i}$ of the components ${C}_{i1}$, \dots ${C}_{is_i}$ intersect one and the same component $E_i$ of the exceptional divisor $D$;
 \item[3)] for $i_1\ne i_2$ the strict transforms $\widetilde{C}_{i_1j_1}$ and $\widetilde{C}_{i_2j_2}$ intersect different components of $D$ (one can say that $E_1$, \dots, $E_r$ are part of the set $\{ E_\sigma \, : \, \sigma \in \Sigma\}$ of irreducible components of $D$).
\end{enumerate}

For $i=1, \ldots, r$, $j=1, \ldots, s_i$, let $\varphi_{ij}:(\CC,0)\to(\CC^2,0)$ be a uniformization of the component $C_{ij}$. For $g\in \calO_{\CC^2,0}$ let $v_{ij}(g)$ be the order of vanishing of $g\circ\varphi_{ij}$ at the origin. The function $v_{ij}(g):\calO_{\CC^2,0}\to\ZZ_{\ge 0}\cup\{+\infty\}$ is a valuation on $\calO_{\CC^2,0}$. Let
$$
v_i(g):=\min\limits_{j=1,\ldots, s_j}v_{ij}(g)\,.
$$
The function $v_i:\calO_{\CC^2,0}\to\ZZ_{\ge 0}\cup\{+\infty\}$ is, in general, not a valuation, but an order function (if $s_i>1$).

The order functions $v_1$, \dots, $v_r$ define in the usual way an $r$-index filtration on $\calO_{\CC^2,0}$:
\begin{equation}\label{curvefilt}
J(\vv)=\{g\in\calO_{\CC^2,0}:\vv(g)\ge\vv\}\,,
\end{equation}
where, as usual, $\vv=(v_1, \ldots, v_r)\in\ZZ_{\ge0}^r$, $\vv(g)=(v_1(g), \dots, v_r(g))$. We shall call it the {\em generalized divisorial filtration}.

Let $\{E_{\sigma}\, : \,  \sigma\in\Sigma\}$ be the set of all irreducible components of the exceptional divisor $D$ ($\Sigma\supset\{1, \ldots, r\}$). Each component $E_{\sigma}$ is isomorphic to the complex projective line $\CC\PP^1$.
For $\sigma\in\Sigma$, let ${\stackrel{\bullet}{E}}_\sigma$ be the ``smooth part" of the
component $E_\sigma$ in the exceptional divisor $D$, that is $E_\sigma$ itself minus the intersection points with
all the other components of $D$, and let ${\stackrel{\circ}{E}}_\sigma$ be the ``smooth part" of the
component $E_\sigma$ in the total transform of the curve $C$, that is $E_\sigma$ itself minus the intersection points with
other components of $D$ and also with the strict transform of the curve $C$. (One has ${\stackrel{\circ}{E}}_\sigma = {\stackrel{\bullet}{E}}_\sigma$ for $\sigma\notin\{1, \ldots, r\}$; for $\sigma=i\in\{1, \ldots, r\}$, ${\stackrel{\circ}{E}}_\sigma$ is ${\stackrel{\bullet}{E}}_\sigma$ minus $s_i$ points.)

For $\sigma\in\Sigma$, let $\widetilde{L}_{\sigma}$ be a smooth arc on the space $X$ of the resolution transversal to $E_{\sigma}$ at a smooth point (i.e. at a point of ${\stackrel{\bullet}{E}}_\sigma$). Let the (irreducible) curve $L_{\sigma}=\pi(\widetilde{L}_{\sigma})$ be given by an equation $g_{\sigma}=0$ ($g_{\sigma}\in\calO_{\CC^2,0}$). The curve $L_{\sigma}$ (or sometimes the function $g_{\sigma}$) is called {\em a curvette} at $E_\sigma$. Let $m_{\sigma\delta}$ ($\sigma, \delta\in\Sigma$) be the order of vanishing of $g_{\sigma}$ along the component $E_\delta$, that is $m_{\sigma\delta}=w_{\delta}(g_{\sigma})$. One can show that $m_{\sigma\delta}=m_{\delta\sigma}$ and the matrix $(m_{\sigma\delta})$ is minus the inverse matrix of the intersection matrix $(E_\sigma\circ E_\delta)$ of the components $E_\sigma$ on the manifold $X$. For $\sigma\in\Sigma$, let $\mm_{\, \sigma}:=(m_{\sigma 1}, \ldots, m_{\sigma r}) \in \ZZ_{\ge0}^r$.

\begin{theorem} \label{theo1}
The Poincar\'e series of the generalized divisorial filtration ({\rm \ref{curvefilt}}) is equal to
\begin{equation}\label{eqntheo}
P_{\{v_i\}}(\ttt)=\prod_{\sigma\in\Sigma}(1-\ttt^{\mm_\sigma})^{-\chi({\stackrel{\bullet}{E}}_\sigma)}\cdot
\prod_{i=1}^r (1-\ttt^{s_i\mm_i})\,.
\end{equation}
\end{theorem}

\begin{example}
Let $s_i=1$ for $i=1, \ldots, r$. In this case $\chi({\stackrel{\circ}{E}}_\sigma)=\chi({\stackrel{\bullet}{E}}_\sigma)$ for $\sigma\notin\{1, \ldots, r\}$ and $\chi({\stackrel{\circ}{E}}_i)=\chi({\stackrel{\bullet}{E}}_i)-1$. Therefore one has
$$
P_{\{v_i\}}(\ttt)=\prod_{\sigma\in\Sigma}(1-\ttt^{\mm_\sigma})^{-\chi({\stackrel{\circ}{E}}_\sigma)}\,.
$$
This is just the formula from \cite{IJM}.
\end{example}

Let $\pi:(X,D)\to(\CC^2,0)$ be a toric resolution corresponding to the Newton diagram $\Gamma=\Gamma_f$ of a ($\Gamma$-non-degenerate) germ $f\in \calO_{\CC^2,0}$. The dual graph of the resolution $\pi$ is a chain. The extreme vertices of this graph correspond to the components of the exceptional divisor intersecting the strict transforms of the coordinate lines in $\CC^2$. (Therefore $\{x=0\}$ and $\{y=0\}$ are curvettes corresponding to these components.) For these two components one has $\chi({\stackrel{\bullet}{E}}_\sigma)=1$, for all others $\chi({\stackrel{\bullet}{E}}_\sigma)=0$. Therefore one has

\begin{corollary}
The Poincar\'e series of  the filtration associated with the Newton diagram $\Gamma$ and defined by the order function $\widehat{v}_i$ corresponding to the facets of $\Gamma$ is equal to
\begin{equation}
P_{\{v_i\}}(\ttt)=\frac{\prod_{i=1}^r (1-\ttt^{s_i\mm_i})}{(1-\ttt^{\vv(x)})(1-\ttt^{\vv(y)})}\,.
\end{equation}
\end{corollary}

\begin{remark}
A function germ $f$ which is non-degenerate with respect to the Newton diagram $\Gamma=\Gamma_f$ can be represented in the form
$f=x^ay^b\prod\limits_{i=1}^rf_i$ where $\{f_i=0\}$ is the union of the components of $\{f=0\}$ whose strict transforms intersect the component $E_i$ of the exceptional divisor of a toric resolution. One can see that the number $s_i$ of irreducible factors in a decomposition of $f_i$ is equal to the integer length of the facet $\gamma_i$ (i.e.\ to the number of integer points in its interior plus one) and
the Newton diagram $\Gamma_i$ of $f_i$ is just the facet $\gamma_i$ of $\Gamma$ translated to the origin inside the positive octant as far as possible. Moreover, the $j$th component of $s_i\mm_i$ is equal to $\min\limits_{\bar{k}\in\Gamma_i}\ell_j(\bar{k})$.
\end{remark}

\begin{proof*}{Proof of Theorem~\ref{theo1}}
Let $\calJ^N_{\CC^2,0}= \calO_{\CC^2,0}/\mathfrak{m}^{N+1}$ be the space of $N$-jets of functions on $(\CC^2,0)$ ($\mathfrak{m}$ is the maximal ideal of $\calO_{\CC^2,0}$). One can see that for a function $g \in \calO_{\CC^2,0}$ with $w_\sigma(g) \leq N$ for all $\sigma \in \Sigma$, the values $w_\sigma(g)$ and also $v_i(g)$ are defined by the $N$-jet $j^Ng$ of $g$. (This follows from the fact that, for $h \in \mathfrak{m}^{N+1}$, all $w_\sigma(h)$ and $v_i(h)$ are greater than $N$.) Let $\widehat{\calJ}^N \subset \calJ^N_{\CC^2,0}$ be the set of $N$-jets $g$ with $w_\sigma(g) \leq N$ for all $\sigma \in \Sigma$. The equation (\ref{Poincare}) implies that 
$$P_{\{ v_i\}}(\ttt) \equiv \int_{\PP \widehat{\calJ}^N} \ttt^{\vv(g)} d \chi$$
modulo terms of degree $>N$. Recall that here $t_i^\infty$ should be assumed to be equal to 0.

Without loss of generality, we can suppose that, for any function $g \in \calO_{\CC^2,0}$ with $w_\sigma(g) \leq N$ for all $\sigma \in \Sigma$, the strict transform $\widetilde{\{ g=0\}}$ of the zero level curve of $g$ intersects the exceptional divisor $D$ only at smooth points, i.e.\ at points of ${\stackrel{\bullet}{D}}= \bigcup_\sigma {\stackrel{\bullet}{E}}_\sigma$. Such a resolution can be obtained, if necessary, by additional blow-ups of intersection points of the components of $D$. Each such blow-up produces an additional component $E_\sigma$ with $\chi({\stackrel{\bullet}{E}}_\sigma)=0$ and therefore it does not effect the right hand side of the equation (\ref{eqntheo}).

Let 
$$
Y= \coprod_{\{ k_\sigma \}} \left( \prod_\sigma S^{k_\sigma} {\stackrel{\bullet}{E}}_\sigma \right) = \prod_\sigma \left( \coprod_{k=0}^\infty S^k{\stackrel{\bullet}{E}}_\sigma \right)
$$
be the configuration space of all effective divisors on ${\stackrel{\bullet}{D}} = \bigcup {\stackrel{\bullet}{E}}_\sigma$ and let $\ww : Y \to \ZZ_{\geq 0}^r$ be the function which maps the component 
$\prod_\sigma S^{k_\sigma} {\stackrel{\bullet}{E}}_\sigma$ of $Y$ to $\sum_\sigma k_\sigma \mm_\sigma$. For a function  $g \in \calO_{\CC^2,0}$ with $w_\sigma(g) \leq N$ for all $\sigma \in \Sigma$, let $I(g) \in Y$ be the intersection of the strict transform $\widetilde{\{ g=0\}}$ of $\{ g=0\}$ with $D$, i.e.\ the collection of the intersection points with multiplicities. One can see that 
$I(g)$ only depends on the $N$-jet of $g$, 
$(w_1(g), \ldots , w_r(g))=\ww(I(g))$ and also $(v_1(g), \ldots , v_r(g))=\ww(I(g))$ if (and only if) for each $i=1, \ldots , r$, the effective divisor $I(g)$ does not contain all the points $p_{i1}, \ldots , p_{is_i}$. (If $I(g)$ contains all the points $p_{i1}, \ldots , p_{is_i}$, then $v_i(g)$ is not determined by $I(g)$.)

For a component $E_\sigma$ of $D$ let $g_{\sigma q}=g_{\sigma q}(x,y)$ be an analytic family of functions  such that $\{ g_{\sigma q} =0\}$ is a curvette corresponding to the component $E_\sigma$ and its strict transform passes through the point $q \in {\stackrel{\bullet}{E}}_\sigma$. (To construct such a family one can take a family of arcs $L_{\sigma q}$ which is analytic in $q \in {\stackrel{\bullet}{E}}_\sigma$ and $L_{\sigma q}$ is transversal to $E_\sigma$ at  the point $q$ and to define $g_{\sigma q}$ as an equation of the corresponding image in $\CC^2 \times {\stackrel{\bullet}{E}}_\sigma$.)

If $A = B \coprod C$, then 
$$
\coprod_{k=0}^\infty S^k A = \left( \coprod_{k=0}^\infty S^k B \right) \times \left(  \coprod_{k=0}^\infty S^k C \right).
$$
This permits to rewrite $Y$ as $Y' \times Y''$, where
$$
Y' =  \prod_\sigma \left( \coprod_{k=0}^\infty S^k{\stackrel{\circ}{E}}_\sigma \right), \quad Y'' = \prod_i \left( \coprod_{k=0}^\infty S^k P_i \right),
$$
where $P_i$ is the set $\{ p_{i1}, \ldots , p_{is_i} \}$ consisting of $s_i$ points. 

For $y \in Y$, $y = \sum_{\sigma, j} \ell_{\sigma j}' q_{\sigma j} + \sum_{i=1}^r \sum_{j=1}^{s_i} \ell_{ij}'' p_{ij}$, where $q_{\sigma j} $ are points of ${\stackrel{\circ}{E}}_\sigma$, let 
$$
g_y := \prod_{\sigma, j} g^{\ell_{\sigma j}'}_{\sigma q_{\sigma j}} \cdot  \prod_{i=1}^r \prod_{j=1}^{s_i} f^{\ell_{ij}''}_{ij},
$$
where $g_{\sigma q_{\sigma j}}$ is the curvette corresponding to $E_\sigma$ through the point $q_{\sigma j}$. One can see that $I(g_y)=y$.

For an element $g \in  \widehat{\calJ}^N$ with $I(g)=y$, one has $I(g)=I(g_y)$, i.e.\ the strict transforms of $\{ g=0 \}$ and $\{ g_y=0 \}$ intersect the exceptional divisor $D$ at the same points with the same multiplicities. This means that the ratio $g_y \circ \pi/g \circ \pi$ of the liftings of $g$ and $g_y$ is regular (has no zeros and poles) on $D$ and therefore it is constant (say, equal to $a$) on it. If $g \neq g_y$, let $h_\lambda:= g_y + \lambda(ag-g_y)$ for $\lambda \in \CC^\ast$. One can see that $w_\sigma(h_\lambda)$ and $v_i(h_\lambda)$ do not depend on $\lambda$. In this way we decompose the space of elements of $\PP \widehat{\calJ}^N$ different from all $g_y$ into $\CC^\ast$-families with constant values of $\vv$. Since the Euler characteristic of $\CC^\ast$ is equal to zero, this means that the integral (with respect to the Euler characteristic) of $\ttt^{\vv}$ over the complement of $\{ g_y \}$ is equal to zero and therefore (up to terms of degree $>N$)
$$
P_{\{ v_i\}}(\ttt) \equiv \int_Y \ttt^{\, \vv(g_y)} d \chi.
$$

For $y\in Y$, $v_i(g_y)$ is finite if and only if $y$ does not contain all the points $p_{i,1}, \dots, p_{i,s_i}$. If, for each $i$, $y$ does not contain all the points $p_{i,1}, \dots, p_{i,s_i}$, one has $\vv(g_y)=\ww(y)$. Therefore
\begin{equation}\label{intprod}
\int\limits_{Y}\ttt^{\, \vv(g_y)}d\chi=
\int\limits_{Y'}\ttt^{\, \ww(y')}d\chi\cdot\int\limits_{Y_0''}\ttt^{\, \ww(y'')}d\chi\,,
\end{equation}
where $Y_0''\subset Y''$ is the set of elements $\sum\limits_{i=1}^r\sum\limits_{j=1}^{s_j}\ell_{ij}p_{ij}$ such that for each $i$ at least one of the coefficients $\ell_{ij}$ is equal to zero.

One has 
$$
\int\limits_{Y'}\ttt^{\, \ww(y')}d\chi=\prod_{\sigma\in\Sigma}\left(\sum_{k=0}^\infty\chi(S^k {\stackrel{\circ}{E}}_\sigma)\ttt^{\, k\mm_\sigma}\right)\,.
$$
Due to the equation
$$
\sum_{k=0}^\infty\chi(S^k Z)t^k=(1-t)^{-\chi(Z)}
$$
one has
\begin{equation}\label{int1}
\int\limits_{Y''}\ttt^{\, \ww(y'')}d\chi=\prod(1-\ttt^{\, \mm_\sigma})^{-\chi({\stackrel{\circ}{E}}_\sigma)}\,.
\end{equation}
(This is just the computation from \cite{IJM}.)

For the second factor in (\ref{intprod}) one has
\begin{eqnarray}\label{int2}
\lefteqn{\int\limits_{Y_0'}
\ttt^{\, \ww(y')}d\chi=\prod_{i=1}^r\left(\sum_{(\ell_{i1},\ldots,\ell_{is_i})\in
\ZZ_{\ge0}^{s_i}\setminus \ZZ_{>0}^{s_i}}\ttt^{\, (\sum\ell_{ij})\mm_i}\right)} \nonumber\\
& = & 
\prod_{i=1}^r
\left(\sum_{(\ell_{i1},\ldots,\ell_{is_i})\in
\ZZ_{\ge0}^{s_i}}\ttt^{\, (\sum\ell_{ij})\mm_i}-
\sum_{(\ell_{i1},\ldots,\ell_{is_i})\in
\ZZ_{>0}^{s_i}}\ttt^{\, (\sum\ell_{ij})\mm_i}
\right) \\
& = & \prod_{i=1}^r\left[(1-\ttt^{\, \mm_i})^{-s_i}- \ttt^{\, s_i\mm_i}(1-\ttt^{\, \mm_i})^{-s_i}\right]=
\prod_{i=1}^r(1-\ttt^{\, \mm_i})^{-s_i}(1- \ttt^{\, s_i\mm_i})\,. \nonumber
\end{eqnarray}
Since $\chi({\stackrel{\bullet}{E}}_\sigma)=\chi({\stackrel{\circ}{E}}_\sigma)+s_i$, the equations (\ref{intprod}), (\ref{int1}), and (\ref{int2}) imply (\ref{eqntheo}).
\end{proof*}

\begin{remark}
Here, in contrast to \cite{IJM}, we make computations of integrals with respect to the Euler characteristic not over $\PP\calO_{\CC^2,0}$, but over a subspace of $\PP{\calJ}^N_{\CC^2,0}$ since the set of function $\{g_y\vert y\in Y\}$ is not measurable in $\PP\calO_{\CC^2,0}$ (i.e.\ its Euler characteristic is not defined).
\end{remark}


\bigskip
\noindent Leibniz Universit\"{a}t Hannover, Institut f\"{u}r Algebraische Geometrie,\\
Postfach 6009, D-30060 Hannover, Germany \\
E-mail: ebeling@math.uni-hannover.de\\

\medskip
\noindent Moscow State University, Faculty of Mechanics and Mathematics,\\
Moscow, GSP-1, 119991, Russia\\
E-mail: sabir@mccme.ru

\end{document}